\documentclass{article}%
\usepackage{graphicx}
\usepackage[intlimits]{amsmath}
\usepackage{latexsym}
\usepackage{amsfonts}
\usepackage{amssymb}%
\newtheorem{theorem}{Theorem}

\newenvironment{proof}[1][Proof]{\noindent{\textbf {#1}  }}  {\hfill$\Box$\bigskip}

\begin{document}

\title{Degree powers in graphs with forbidden subgraphs}
\author{B\'{e}la Bollob\'{a}s\thanks{Department of Mathematical Sciences, University
of Memphis, Memphis TN 38152, USA} \thanks{Trinity College, Cambridge CB2 1TQ,
UK} \thanks{Research supported in part by DARPA grant F33615-01-C-1900.} \ and
Vladimir Nikiforov$^{\ast}$}
\maketitle

\begin{abstract}
For every real $p>0$ and simple graph $G,$ set
\[
f\left(  p,G\right)  =\sum_{u\in V\left(  G\right)  }d^{p}\left(  u\right)  ,
\]
and let $\phi\left(  p,n,r\right)  $ be the maximum of $f\left(  p,G\right)  $
taken over all $K_{r+1}$-free graphs $G$ of order $n.$ We prove that, if
$0<p<r,$ then%
\[
\phi\left(  p,n,r\right)  =f\left(  p,T_{r}\left(  n\right)  \right)  ,
\]
where $T_{r}\left(  n\right)  $ is the $r$-partite Turan graph of order $n$.
For every $p\geq r+\left\lceil \sqrt{2r}\right\rceil $ and $n$ large, we show
that%
\[
\phi\left(  p,n,r\right)  >\left(  1+\varepsilon\right)  f\left(
p,T_{r}\left(  n\right)  \right)
\]
for some $\varepsilon=\varepsilon\left(  r\right)  >0.$

Our results settle two conjectures of Caro and Yuster.

\end{abstract}

\section{Introduction}

Our notation and terminology are standard (see, e.g. \cite{Bol98}).

Caro and Yuster \cite{CaYu00} introduced and investigated the function
\[
f\left(  p,G\right)  =\sum_{u\in V\left(  G\right)  }d^{p}\left(  u\right)  ,
\]
where $p\geq1$ is integer and $G$ is a graph. Writing $\phi\left(
r,p,n\right)  $ for the maximum value of $f\left(  p,G\right)  $ taken over
all $K_{r+1}$-free graphs $G$ of order $n,$ Caro and Yuster stated that, for
every $p\geq1,$
\begin{equation}
\phi\left(  r,p,n\right)  =f\left(  p,T_{r}\left(  n\right)  \right)  ,
\label{YC}%
\end{equation}
where $T_{r}\left(  n\right)  $ is the $r$-partite Tur\'{a}n graph of order
$n$. However, simple examples show that (\ref{YC}) fails for every fixed
$r\geq2$ and all sufficiently large $p$ and $n;$ this was observed by Schelp
\cite{Sch01}. A natural problem arises: given $r\geq2,$ determine those real
values $p>0,$ for which equality (\ref{YC}) holds. Furthermore, determine the
asymptotic value of $\phi\left(  r,p,n\right)  $ for large $n.$

In this note we essentially answer these questions. In Section \ref{S1} we
prove that (\ref{YC}) holds whenever $0<p<r$ and $n$ is large. Next, in
Section \ref{S2}, we describe the asymptotic structure of $K_{r+1}$-free
graphs $G$ of order $n$ such that $f\left(  p,G\right)  =\phi\left(
r,p,n\right)  .$ We deduce that, if $p\geq r+\left\lceil \sqrt{2r}\right\rceil
$ and $n$ is large, then%
\[
\phi\left(  r,p,n\right)  >\left(  1+\varepsilon\right)  f\left(
p,T_{r}\left(  n\right)  \right)
\]
for some $\varepsilon=\varepsilon\left(  r\right)  >0.$ This disproves
Conjecture 6.2 in \cite{CaYu00}. In particular,
\[
\frac{r}{pe}\geq\frac{\phi\left(  r,p,n\right)  }{n^{p+1}}\geq\frac
{r-1}{\left(  p+1\right)  e}%
\]
holds for large $n,$ and therefore, for any fixed $r\geq2,$
\[
\lim_{n\rightarrow\infty}\frac{\phi\left(  r,p,n\right)  }{f\left(
p,T_{r}\left(  n\right)  \right)  }%
\]
grows exponentially in $p.$

The case $r=2$ is considered in detail in Section \ref{S3}; we show that, if
$r=2,$ equality (\ref{YC}) holds for $0<p\leq3,$ and is false for every $p>3$
and $n$ large.

In Section \ref{S4} we extend the above setup. For a fixed $\left(
r+1\right)  $-chromatic graph $H,$ $\left(  r\geq2\right)  ,$ let $\phi\left(
H,p,n\right)  $ be the maximum value of $f\left(  p,G\right)  $ taken over all
$H$-free graphs $G$ of order $n.$ It turns out that, for every $r$ and $p,$
\begin{equation}
\phi\left(  H,p,n\right)  =\phi\left(  r,p,n\right)  +o\left(  n^{p+1}\right)
. \label{YC1}%
\end{equation}
This result completely settles, with the proper changes, Conjecture 6.1 of
\cite{CaYu00}. In fact, Pikhurko \cite{Pik01} proved this for $p\geq1,$
although he incorrectly assumed that (\ref{YC}) holds for all sufficiently
large $n.$

\section{\label{S1}The function $\phi\left(  r,p,n\right)  $ for $p<r$}

In this section we shall prove the following theorem.

\begin{theorem}
\label{mainT}For every $r\geq2,$ $0<p<r,$ and sufficiently large $n$,
\[
\phi\left(  r,p,n\right)  =f\left(  p,T_{r}\left(  n\right)  \right)  .
\]

\end{theorem}

\begin{proof}
Erd\H{o}s \cite{Erd70} proved that, for every $K_{r+1}$-free graph $G,$ there
exists an $r$-partite graph $H$ with $V\left(  H\right)  =V\left(  G\right)  $
such that $d_{G}\left(  u\right)  \leq d_{H}\left(  u\right)  $ for every
$u\in V\left(  G\right)  $. As Caro and Yuster noticed, this implies that, for
$K_{r+1}$-free graphs $G$ of order $n,$ if $f\left(  p,G\right)  $ attains a
maximum then $G$ is a complete $r$-partite graph. Every complete $r$-partite
graph is defined uniquely by the size of its vertex classes, that is, by a
vector $\left(  n_{i}\right)  _{1}^{r}$ of positive integers satisfying
$n_{1}+...+n_{r}=n;$ note that the Tur\'{a}n graph $T_{r}\left(  n\right)  $
is uniquely characterized by the condition $\left\vert n_{i}-n_{j}\right\vert
\leq1$ for every $i,j\in\left[  r\right]  .$ Thus we have\
\begin{equation}
\phi\left(  r,p,n\right)  =\max\left\{  \sum_{i=1}^{r}n_{i}\left(
n-n_{i}\right)  ^{p}\text{ }:\text{\ }n_{1}+...+n_{r}=n,\text{ }1\leq
n_{1}\leq...\leq n_{r}\right\}  . \label{bas}%
\end{equation}

Let $\left(  n_{i}\right)  _{1}^{r}$ be a vector on which the value of
$\phi\left(  r,p,n\right)  $ is attained. Routine calculations show that the
function $x\left(  n-x\right)  ^{p}$ increases for $0\leq x\leq\frac{n}{p+1},$
decreases for $\frac{n}{p+1}\leq x\leq n,$ and is concave for $\frac{2n}%
{p+1}\leq x\leq n.$ If $n_{r}\leq\left\lfloor \frac{2n}{p+1}\right\rfloor $,
the concavity of $x\left(  n-x\right)  ^{p}$ implies that $n_{1}-n_{r}\leq1,$
and the proof is completed, so we shall assume $n_{r}>\left\lfloor \frac
{2n}{p+1}\right\rfloor .$ Hence we deduce
\begin{equation}
n_{1}\left(  r-1\right)  +\left\lfloor \frac{2n}{p+1}\right\rfloor <\text{
}n_{1}+...+n_{r}=n. \label{in1}%
\end{equation}
We shall also assume
\begin{equation}
n_{1}\geq\left\lfloor \frac{n}{p+1}\right\rfloor , \label{in2}%
\end{equation}
since otherwise, adding $1$ to $n_{r}$ and subtracting $1$ from $n_{1},$ the
value $\sum_{i=1}^{r}n_{i}\left(  n-n_{i}\right)  ^{p}$ will increase,
contradicting the choice of $\left(  n_{i}\right)  _{1}^{r}$. Notice that, as
$n_{1}\leq n/r,$ inequality (\ref{in2}) is enough to prove the assertion for
$p\leq r-1$ and every $n.$ From (\ref{in1}) and (\ref{in2}), we obtain that%
\[
\left(  r-1\right)  \left\lfloor \frac{n}{p+1}\right\rfloor +\left\lfloor
\frac{2n}{p+1}\right\rfloor <n.
\]
Letting $n\rightarrow\infty,$ we see that $p\geq r,$ contradicting the
assumption and completing the proof.
\end{proof}

Maximizing independently each summand in (\ref{bas}), we see that, for every
$r\geq2$ and $p>0$,%
\begin{equation}
\phi\left(  r,p,n\right)  \leq\frac{r}{p+1}\left(  \frac{p}{p+1}\right)
^{p}n^{p+1}. \label{filo}%
\end{equation}

\section{\label{S2}The asymptotics of $\phi\left(  r,p,n\right)  $}

In this section we find the asymptotic structure of $K_{r+1}$-free graphs $G$
of order $n$ satisfying $f\left(  p,G\right)  =\phi\left(  r,p,n\right)  ,$
and deduce asymptotic bounds on $\phi\left(  r,p,n\right)  .$

\begin{theorem}
For all $r\geq2$ and $p>0,$ there exists $c=c\left(  p,r\right)  $ such that
the following assertion holds.

If $f\left(  p,G\right)  =\phi\left(  r,p,n\right)  $ for some $K_{r+1}$-free
graph $G$ of order $n,$ then $G$ is a complete $r$-partite graph having $r-1$
vertex classes of size $cn+o\left(  n\right)  .$
\end{theorem}

\begin{proof}
We already know that $G$ is a complete $r$-partite graph; let $n_{1}%
\leq...\leq n_{r}$ be the sizes of its vertex classes and, for every
$i\in\left[  r\right]  ,$ set $y_{i}=n_{i}/n$. It is easy to see that%
\[
\phi\left(  r,p,n\right)  =\psi\left(  r,p\right)  n^{p+1}+o\left(
n^{p+1}\right)  ,
\]
where the function $\psi\left(  r,p\right)  $ is defined as
\[
\psi\left(  r,p\right)  =\max\left\{  \sum_{i=1}^{r}x_{i}\left(
1-x_{i}\right)  ^{p}\text{ }:\text{\ }x_{1}+...+x_{r}=1,\text{ }0\leq
x_{1}\leq...\leq x_{r}\right\}
\]

We shall show that if the above maximum is attained at $\left(  x_{i}\right)
_{1}^{r},$ then $x_{1}=...=x_{r-1}.$ Indeed, the function $x\left(
1-x\right)  ^{p}$ is concave for $0\leq x\leq2/\left(  p+1\right)  ,$ and
convex for $2/\left(  p+1\right)  \leq x\leq1.$ Hence, there is at most one
$x_{i}$ in the interval $(2/\left(  p+1\right)  \leq x\leq1],$ which can only
be $x_{r}.$ Thus $x_{1},...,x_{r-1}$ are all in the interval $\left[
0,2/\left(  p+1\right)  \right]  ,$ and so, by the concavity of $x\left(
1-x\right)  ^{p}$, they are equal. We conclude that, if
\[
0\leq x_{1}\leq...\leq x_{r},\text{ }x_{1}+...+x_{r}=1,
\]
and $x_{j}>x_{i}$ for some $1\leq i<j\leq r-1$, then $\sum_{i=1}^{r}%
x_{i}\left(  1-x_{i}\right)  ^{p}$ is below its maximum value. Applying this
conclusion to the numbers $\left(  y_{i}\right)  _{1}^{r},$ we deduce the
assertion of the theorem.
\end{proof}

Set
\[
g\left(  r,p,x\right)  =\left(  r-1\right)  x\left(  1-x\right)  ^{p}+\left(
1-\left(  r-1\right)  x\right)  \left(  rx\right)  ^{p}.
\]
From the previous theorem it follows that
\[
\psi\left(  r,p\right)  =\max_{0\leq x\leq1/\left(  r-1\right)  }g\left(
r,p,x\right)  .
\]
Finding $\psi\left(  r,p\right)  $ is not easy when $p>r.$ In fact, for some
$p>r,$ there exist $0<x<y<1$ such that%
\[
\psi\left(  r,p\right)  =g\left(  r,p,x\right)  =g\left(  r,p,y\right)  .
\]

In view of the original claim concerning (\ref{YC}), it is somewhat
surprising, that for $p>2r-1$, the point $x=1/r,$ corresponding to the
Tur\'{a}n graph, not only fails to be a maximum of $g\left(  r,p,x\right)  $,
but, in fact, is a local minimum.

Observe that
\[
\frac{f\left(  p,T_{r}\left(  n\right)  \right)  }{n^{p+1}}=\left(  \frac
{r-1}{r}\right)  ^{p}+o\left(  1\right)  ,
\]
so, to find for which $p$ the function $\phi\left(  r,p,n\right)  $ is
significantly greater than $f\left(  p,T_{r}\left(  n\right)  \right)  $, we
shall compare $\psi\left(  r,p\right)  $ to $\left(  \frac{r-1}{r}\right)
^{p}$.

\begin{theorem}
\label{th2}Let $r\geq2,$ $p\geq r+\left\lceil \sqrt{2r}\right\rceil .$ Then%
\[
\psi\left(  r,p\right)  >\left(  1+\varepsilon\right)  \left(  \frac{r-1}%
{r}\right)  ^{p}.
\]
for some $\varepsilon=\varepsilon\left(  r\right)  >0.$
\end{theorem}

\begin{proof}
We have
\begin{align*}
\psi\left(  r,p\right)   &  \geq g\left(  r,p,\frac{1}{p}\right)  =\frac
{r-1}{p}\left(  \frac{p-1}{p}\right)  ^{p}+\left(  1-\frac{r-1}{p}\right)
\left(  \frac{r-1}{p}\right)  ^{p}\\
&  >\frac{r-1}{p}\left(  \frac{p-1}{p}\right)  ^{p}.
\end{align*}

To prove the theorem, it suffices to show that
\begin{equation}
\frac{r-1}{p}\left(  \frac{\left(  p-1\right)  r}{p\left(  r-1\right)
}\right)  ^{p}>1+\varepsilon\label{in3}%
\end{equation}
for some $\varepsilon=\varepsilon\left(  r\right)  >0.$ Routine calculations
show that
\[
\frac{r-1}{p}\left(  1+\frac{p-r}{p\left(  r-1\right)  }\right)  ^{p}%
\]
increases with $p.$ Thus, setting $q=\left\lceil \sqrt{2r}\right\rceil ,$ we
find that%
\begin{align*}
&  \frac{r-1}{p}\left(  1+\frac{p-r}{p\left(  r-1\right)  }\right)  ^{p}\\
&  \geq\frac{r-1}{r+q}\left(  1+\binom{r+q}{1}\frac{q}{\left(  r+q\right)
\left(  r-1\right)  }+\binom{\left(  r+q\right)  }{2}\frac{q^{2}}{\left(
r+q\right)  ^{2}\left(  r-1\right)  ^{2}}\right) \\
&  =\frac{r-1}{r+q}+\frac{q}{r+q}+\frac{q^{2}\left(  r+q-1\right)  }{2\left(
r+q\right)  ^{2}\left(  r-1\right)  }\geq1-\frac{1}{r+q}+\frac{r\left(
r+q-1\right)  }{\left(  r+q\right)  ^{2}\left(  r-1\right)  }\\
&  =1+\frac{r\left(  r+q-1\right)  -\left(  r+q\right)  \left(  r-1\right)
}{\left(  r+q\right)  ^{2}\left(  r-1\right)  }=1+\frac{q}{\left(  r+q\right)
^{2}\left(  r-1\right)  }.
\end{align*}
Hence, (\ref{in3}) holds with%
\[
\varepsilon=\frac{\left\lceil \sqrt{2r}\right\rceil }{\left(  r+\left\lceil
\sqrt{2r}\right\rceil \right)  ^{2}\left(  r-1\right)  },
\]
completing the proof.
\end{proof}

We have, for $n$ sufficiently large,
\begin{align*}
\frac{\phi\left(  r,p,n\right)  }{n^{p+1}}  &  =\psi\left(  r,p\right)
+o\left(  1\right)  \geq g\left(  r,p,\frac{1}{p+1}\right)  +o\left(  1\right)
\\
&  =\frac{r-1}{p+1}\left(  \frac{p}{p+1}\right)  ^{p}+\left(  1-\frac
{r-1}{p+1}\right)  \left(  \frac{r-1}{p+1}\right)  ^{p}+o\left(  1\right) \\
&  >\frac{r-1}{p+1}\left(  \frac{p}{p+1}\right)  ^{p}.
\end{align*}
Hence, in view of (\ref{filo}), we find that, for $n$ large,
\[
\frac{r}{pe}\geq\frac{r}{p}\left(  \frac{p}{p+1}\right)  ^{p+1}\geq\frac
{\phi\left(  r,p,n\right)  }{n^{p+1}}\geq\frac{r-1}{p+1}\left(  \frac{p}%
{p+1}\right)  ^{p}\geq\frac{\left(  r-1\right)  }{\left(  p+1\right)  e}.
\]

In particular, we deduce that, for any fixed $r\geq2,$
\[
\lim_{n\rightarrow\infty}\frac{\phi\left(  r,p,n\right)  }{f\left(
p,T_{r}\left(  n\right)  \right)  }%
\]
grows exponentially in $p.$

\section{\label{S3}Triangle-free graphs}

For triangle-free graphs, i.e., $r=2$, we are able to pinpoint the value of
$p$ for which (\ref{YC}) fails, as stated in the following theorem.

\begin{theorem}
\label{CY3}If $0<p\leq3$ then%
\begin{equation}
\phi\left(  3,p,n\right)  =f\left(  p,T_{2}\left(  n\right)  \right)
.\label{treq}%
\end{equation}
For every $\varepsilon>0,$ there exists $\delta$ such that if $p>3+\delta$
then
\begin{equation}
\phi\left(  3,p,n\right)  >\left(  1+\varepsilon\right)  f\left(
p,T_{2}\left(  n\right)  \right)  \label{trin}%
\end{equation}
for $n$ sufficiently large.
\end{theorem}

\begin{proof}
We start by proving (\ref{treq}). From the proof of Theorem \ref{mainT} we
know that
\[
\phi\left(  p,n,3\right)  =\max_{k\in\left\lceil n/2\right\rceil }\text{
}\left\{  k\left(  n-k\right)  ^{p}+\left(  n-k\right)  k^{p}\right\}  .
\]

Our goal is to prove that the above maximum is attained at $k=\left\lceil
n/2\right\rceil .$

If $0<p\leq2,$ the function $x\left(  1-x\right)  ^{p}$ is concave, and
(\ref{treq}) follows immediately.

Next, assume that $2<p\leq3;$ we claim that the function
\[
g\left(  x\right)  =\left(  1+x\right)  \left(  1-x\right)  ^{p}+\left(
1-x\right)  \left(  1+x\right)  ^{p}%
\]
is concave for $\left\vert x\right\vert \leq1.$ Indeed, we have
\begin{align*}
g\left(  x\right)   &  =\left(  1-x^{2}\right)  \left(  \left(  1-x\right)
^{p-1}+\left(  1+x\right)  ^{p-1}\right)  =2\left(  1-x^{2}\right)  \sum
_{n=0}^{\infty}\binom{p-1}{2n}x^{2n}\\
&  =2+2\sum_{n=1}^{\infty}\left(  \binom{p-1}{2n}-\binom{p-1}{2n-2}\right)
x^{2n}\\
&  =2+2\sum_{n=1}^{\infty}\binom{p-1}{2n-2}\left(  \frac{\left(
p-2n-1\right)  \left(  p-2n-2\right)  }{\left(  2n-1\right)  2n}-1\right)
x^{2n}.
\end{align*}
Since, for every $n,$ the coefficient of $x^{2n}$ is nonpositive, the function
$g\left(  x\right)  $ is concave, as claimed.

Therefore, the function $h\left(  x\right)  =x\left(  n-x\right)  ^{p}+\left(
n-x\right)  x^{p}$ is concave for $1\leq x\leq n.$ Hence, for every integer
$k\in\left[  n\right]  ,$ we have
\begin{align*}
h\left(  \left\lceil \frac{n}{2}\right\rceil \right)  +h\left(  \left\lfloor
\frac{n}{2}\right\rfloor \right)   &  \geq h\left(  k\right)  +h\left(
n-k\right)  =2h\left(  k\right)  \\
&  =2\left(  k\left(  n-k\right)  ^{p}+\left(  n-k\right)  k^{p}\right)  ,
\end{align*}
proving (\ref{treq}).

Inequality (\ref{trin}) follows easily, since, in fact, for every $p>3,$ the
function $g\left(  x\right)  $ has a local minimum at $0.$
\end{proof}

\section{\label{S4}$H$-free graphs}

In this section we are going to prove the following theorem.

\begin{theorem}
\label{PKT}For every $r\geq2,$ and $p>0,$
\[
\phi\left(  H,p,n\right)  =\phi\left(  r,p,n\right)  +o\left(  n^{p+1}\right)
.
\]

\end{theorem}

A few words about this theorem seem in place. As already noted, Pikhurko
\cite{Pik01} proved the assertion for $p\geq1;$ although he incorrectly
assumed that (\ref{YC}) holds for all $p$ and sufficiently large $n,$ his
proof is valid, since it is independent of the exact value of $\phi\left(
r,p,n\right)  .$ Our proof is close to Pikhurko's, and is given only for the
sake of completeness.

We shall need the following theorem (for a proof see, e.g., \cite{Bol98},
Theorem 33, p. 132).

\begin{theorem}
\label{HT}Suppose $H$ is an $\left(  r+1\right)  $-chromatic graph. Every
$H$-free graph $G$ of sufficiently large order $n$ can be made $K_{r+1}$-free
by removing $o\left(  n^{2}\right)  $ edges.
\end{theorem}

\begin{proof}
[Proof of Theorem \ref{PKT}]Select a $K_{r+1}$-free graph $G$ of order $n$
such that $f\left(  p,G\right)  =\phi\left(  r,p,n\right)  .$ Since $G$ is
$r$-partite, it is $H$-free, so we have $\phi\left(  H,p,n\right)  \geq
\phi\left(  r,p,n\right)  .$ Let now $G$ be a $H$-free graph of order $n$ such
that%
\[
f\left(  p,G\right)  =\phi\left(  H,p,n\right)  .
\]

Theorem \ref{HT} implies that there exists a $K_{r+1}$-free graph $F$ that may
be obtained from $G$ by removing at most $o\left(  n^{2}\right)  $ edges.
Obviously, we have
\[
e\left(  G\right)  =e\left(  F\right)  +o\left(  n^{2}\right)  \leq\frac
{r-1}{2r}n^{2}+o\left(  n^{2}\right)  .
\]

For $0<p\leq1,$ by Jensen's inequality, we have
\[
\left(  \frac{1}{n}f\left(  p,G\right)  \right)  ^{1/p}\leq\frac{1}{n}f\left(
1,G\right)  =\frac{1}{n}2e\left(  G\right)  \leq\frac{r-1}{r}n+o\left(
n\right)  .
\]
Hence, we find that
\[
f\left(  p,G\right)  \leq\left(  \frac{r-1}{r}\right)  ^{p}n^{p+1}+o\left(
n^{p+1}\right)  =\phi\left(  r,p,n\right)  +o\left(  n^{p+1}\right)  ,
\]
completing the proof.

Next, assume that $p>1.$ Since the function $xn^{p-1}-x^{p}$ is decreasing for
$0\leq x\leq n,$ we find that
\[
d_{G}^{p}\left(  u\right)  -d_{F}^{p}\left(  u\right)  \leq\left(
d_{G}\left(  u\right)  -d_{F}\left(  u\right)  \right)  n^{p-1}%
\]
for every $u\in V\left(  G\right)  .$ Summing this inequality for all $u\in
V\left(  G\right)  $, we obtain%
\begin{align*}
f\left(  p,G\right)   &  \leq f\left(  p,F\right)  +\left(  d_{G}\left(
u\right)  -d_{F}\left(  u\right)  \right)  n^{p-1}=f\left(  p,F\right)
+o\left(  n^{p+1}\right) \\
&  \leq\phi\left(  r,p,n\right)  +o\left(  n^{p+1}\right)  ,
\end{align*}
completing the proof.
\end{proof}

\section{Concluding remarks}

It seems interesting to find, for each $r\geq3,$ the minimum $p$ for which the
equality (\ref{YC}) is essentially false for $n$ large. Computer calculations
show that this value is roughly $4.9$ for $r=3,$ and $6.2$ for $r=4$,
suggesting that the answer might not be easy.

\end{document}